# A secondary Chern-Euler class

By Ji-Ping Sha

### Introduction

Let $\xi$ be a smooth oriented vector bundle, with $n$-dimensional fibre, over a smooth manifold $M$. Denote by $\hat{\xi}$ the fibrewise one-point compactification of $\xi$. The main purpose of this paper is to define geometrically a canonical element $\Upsilon(\xi)$ in $H^n(\hat{\xi}, \mathbb{Q})$ ($H^n(\hat{\xi}, \mathbb{Z}) \otimes \frac{1}{2}$, to be more precise). The element $\Upsilon(\xi)$ is a secondary characteristic class to the Euler class in the fashion of Chern-Simons. Two properties of this element are described as follows.

The first one is in a very classical setting. Suppose $\xi$ is the tangent bundle $TM$ of $M$ (hence $M$ is oriented). In this case we denote $\hat{\xi}$ by $\Sigma M$ and simply write $\Upsilon$ for $\Upsilon(\hat{\xi})$.

Suppose $M$ is the boundary of a compact $(n+1)$-dimensional smooth manifold $X$. Let $V$ be a nowhere zero smooth vector field given on $M$ which is tangent to $X$, but not necessarily tangent or transversal to $M$. The vector field $V$ naturally defines a cross section $\alpha : M \to \Sigma M$. One can extend $V$ to a smooth tangent vector field $\overline{V}$ on $X$ with only isolated (hence only a finite number of) zeros. Since such extensions are generic we shall, for convenience, call any such extension a generic extension. At an isolated zero point $p$ of $\overline{V}$, let $\operatorname{ind}_p(\overline{V})$ be the index of $\overline{V}$ at $p$ defined as usual. We then have the following:

THEOREM 0.1. *For any generic extension $\overline{V}$ of $V$, if $p_1, \ldots, p_k$ are the zero points of $\overline{V}$ then*

$$\sum_{j=1}^{k} \operatorname{ind}_{p_j}(\overline{V}) = \begin{cases} \chi(X) + \alpha^*(\Upsilon)([M]) & \text{if } n \text{ is odd} \\ \alpha^*(\Upsilon)([M]) & \text{if } n \text{ is even} \end{cases}$$

*where $\chi(X)$ is the Euler characteristic of $X$.*

Notice that, in case $M$ is empty, if we establish as a convention that $\alpha^*(\Upsilon)([M]) = 0$, then the theorem above is a generalization of a well-known theorem of Poincaré-Hopf (cf. [M]). In general if $M$ is not empty, it is easy to see from the Poincaré-Hopf theorem that the sum $\sum_{j=1}^{k} \operatorname{ind}_{p_j}(\overline{V})$ does not depend on the extension $\overline{V}$; and in case $n$ is even, it does not depend on $X$.

Our theorem above relates the sum to a specific topological invariant of the boundary.

*Note.* Generalizing the Poincaré-Hopf index theorem for vector fields to manifolds with boundary has been studied by C. Pugh and D. Gottlieb (cf. [G], [P]). The formulae obtained in [G] and [P] however do not seem to link directly to the global topological invariant of the boundary in general.

The second property of $\Upsilon(\xi)$ is that it is closely related to the Thom class. Let $\xi_\infty$ be the $\infty$-section of $\hat{\xi}$, and let $\gamma(\xi) \in H^n(\hat{\xi}, \xi_\infty)$, with integer coefficients, be the Thom class of $\xi$. We shall show the following:

THEOREM 0.2. *The natural homomorphism $\jmath^* : H^n(\hat{\xi}, \xi_\infty) \to H^n(\hat{\xi})$ is injective, and*

$$\jmath^*(\gamma(\xi)) = \Upsilon(\xi) + \frac{1}{2}\sigma^*(e(\xi))$$

*where $e(\xi)$ is the Euler class of $\xi$, and $\sigma : \hat{\xi} \to M$ is the projection.*

The construction of $\Upsilon(\xi)$ is explicit, and is inspired by Chern's well-known proof of the Gauss-Bonnet theorem. While $\Upsilon(\xi)$ can be defined formally in a pretty straightforward way, in order to see its nature as a secondary characteristic class and prove Theorem 0.1 above, we shall first construct it as an element in $H^n(\hat{\xi}, \mathbb{R})$ in Section 1; the construction depends on choice of a connection on $\xi$. A proof of Theorem 0.1 is given in Section 2, while the proof of the topological invariance of the $\Upsilon(\xi)$ constructed in Section 1 is postponed to Section 3. There we shall see that $\Upsilon(\xi)$ is defined in $H^n(\hat{\xi}, \mathbb{Z}) \otimes \frac{1}{2}$, and prove Theorem 0.2.

*Acknowledgment.* The author would like to thank the referee for the suggestions that improved the exposition.

## Section 1

In this section we first construct, in a natural way, a closed differential $n$-form $\Psi$ on $\hat{\xi}$ (note that $\hat{\xi}$ has a canonical smooth structure). The form $\Psi$ then represents an element in the de Rham cohomology $H^n(\hat{\xi}, \mathbb{R})$. It will be seen in subsequent sections that this element is in fact half integral, and does not depend on various choices involved in the construction.

The construction of $\Psi$ follows the well-known work of Chern in [C], with some modifications particularly in the case when the dimension $n$ of the fibre of $\xi$ is even. For completeness we shall show the construction in detail, while leaving some needed fundamental background in differential geometry to the references (e.g. [KN]).

To start with, we fix an $SO(n)$-connection $\omega$ on $\xi$, and let $\Omega$ be the curvature. Let us first explain some notational conventions that we are going to use, most of them standard.

We denote by $\langle \, , \, \rangle$ and $\| \, \|$ the underlying metric and the induced norm, respectively, on $\xi$. The same notation will be used for the induced metric and norm on any other vector bundle associated to $\xi$.

Let $\nu$ be the *canonical* trivial oriented real line bundle over $M$ with the trivial connection. Let $E = \nu \oplus \xi$. We then have an obvious (orientation-preserving) diffeomorphism
$$\hat{\xi} \approx \{v \in E : \|v\| = 1\}$$
in which the 0-section of $\hat{\xi}$ is identified with $1 \oplus 0$, the $\infty$-section of $\hat{\xi}$ is identified with $-1 \oplus 0$, and the unit sphere bundle of $\xi$ is in $0 \oplus \xi$. We shall always use this diffeomorphism without further notice.

The obviously induced $SO(n+1)$-connection and curvature on $E$ will still be denoted by $\omega$ and $\Omega$ respectively. Throughout our calculation, we shall choose an oriented local orthonormal frame field for $\xi$ on $M$. Together with the canonical (positive) unit vector of $\nu$ in the first position, this forms the oriented local orthonormal frame field we shall choose for $E$ on $M$. To simplify the notation without causing any ambiguity, we shall view $\omega$ ($\Omega$, resp.) as an $so(n+1)$-valued 1-form (2-form, resp.) on $M$, with respect to the chosen frame field. Recall $\Omega = d\omega + \omega \wedge \omega$, where matrix multiplication is understood. Also notice that the first row and column of $\omega$ and $\Omega$ are always 0.

As in the introduction, let $\sigma : \hat{\xi} \to M$ be the projection. For any differential form $A$ on $M$, for the sake of simplicity, we shall write $A$ for $\sigma^*(A)$ on $\hat{\xi}$ wherever it can be easily understood from the context.

Let $u = \begin{pmatrix} u_1 \\ \vdots \\ u_{n+1} \end{pmatrix}$ be the $\mathbb{R}^{n+1}$-valued function on $\hat{\xi}$, associated to a chosen local frame field $e = (e_1, \ldots, e_{n+1})$ for $E$ described above, defined by
$$v = \sum_{i=1}^{n+1} u_i(v) e_i, \quad \forall v \in \hat{\xi},$$
and let $\theta = \begin{pmatrix} \theta_1 \\ \vdots \\ \theta_{n+1} \end{pmatrix}$ be the $\mathbb{R}^{n+1}$-valued 1-form defined by
$$\theta = du + \omega u.$$

The definition of $u$ and $\theta$ depends on the choice of the local frame field of course. However, if the local frame field $e$ is replaced by any other frame field $eg$ for some $SO(n+1)$-valued local function $g$, then it is easily seen that $u$ and $\theta$ are replaced by $g^{-1}u$ and $g^{-1}\theta$ correspondingly.

We are now ready to define the form $\Psi$. Suppose $n = 2m$ or $2m+1$. Set

$$\Psi_j = \sum_\tau (-1)^\tau u_{\tau(1)} \theta_{\tau(2)} \wedge \cdots \wedge \theta_{\tau(n-2j+1)} \wedge \Omega_{\tau(n-2j+2)\tau(n-2j+3)} \wedge \cdots \wedge \Omega_{\tau(n)\tau(n+1)}$$

for $j = 0, 1, \ldots, m$, where the summation is over all the permutations $\tau$ of $\{1, \ldots, n+1\}$, and $\Omega_{st}$ denotes the $(s,t)$-entry of the matrix $\Omega$ as usual.

It is easy to see that the definition of each of the $\Psi_j$ above does not depend on the choice of local frame, and hence is a globally well defined $n$-form on $\hat\xi$. We now define

$$\Psi = \frac{1}{(n-1)!! c_n} \sum_{j=0}^m \frac{1}{2^j j! (n-2j)!!} \Psi_j$$

where

$$c_n = \begin{cases} \frac{2(2\pi)^m}{(n-1)!!} & \text{if } n = 2m \\ \frac{(2\pi)^{m+1}}{(n-1)!!} & \text{if } n = 2m+1 \end{cases}$$

is the volume of the Euclidean $n$-dimensional sphere $S^n$.

We summarize some basic properties of $\Psi$ in the following proposition. Its proof follows from the computations in [C], and hence is omitted. We state this proposition in the more general setting where $E$ is an arbitrary oriented vector bundle over $M$, with $(n+1)$-dimensional fiber, and $\omega$ is an arbitrary $SO(n+1)$-connection on $E$.

PROPOSITION 1.1.
(1)
$$d\Psi = \begin{cases} 0 & \text{if } n = 2m \\ -E(\Omega) & \text{if } n = 2m+1 \end{cases}$$

where, for $n = 2m+1$,

$$E(\Omega) = \frac{1}{(4\pi)^{m+1}(m+1)!} \sum_\tau (-1)^\tau \Omega_{\tau(1)\tau(2)} \cdots \Omega_{\tau(n)\tau(n+1)}$$

is the Euler curvature form of $E$.

(2) If $\imath : S^n \to \hat\xi$ is any (orientation-preserving) isometry from the euclidean sphere $S^n$ to a fibre of $\sigma : \hat\xi \to M$, then $\imath^*(\Psi) = \frac{1}{c_n}\mathrm{vol}$, where vol denotes the volume form on $S^n$.

Returning to the special case when $E = \nu \oplus \xi$ and $\omega$ is induced from a connection on $\xi$, we have that $\Psi$ is a closed $n$-form on $\hat\xi$, since the first row and column of $\Omega$ are 0.

Finally we note that the construction of $\Psi$ is obviously natural (in the category of oriented vector bundles with Riemannian connection).

## Section 2

In this section we assume $\xi$ is the tangent bundle $TM$ of $M$. Let $\Upsilon$ be the cohomology class in $H^n(\Sigma M, \mathbb{R})$ represented by the $n$-form $\Psi$ constructed in last section. We now prove Theorem 0.1 stated in the introduction. First we note the following:

*Remark* 2.1. The vector bundle $\nu \oplus TM$ can naturally be viewed as one over $\mathbb{R} \times M$, and identified with the tangent bundle $T(\mathbb{R} \times M)$. The $SO(n+1)$-connection $\omega$ in Section 1 is then associated with the Riemannian product metric on $\mathbb{R} \times M$.

Suppose $M$ is the boundary of a compact $(n+1)$-dimensional manifold $X$. Assume $X$ is orientable. We orient $X$ consistently with the orientation of $M$. By Remark 2.1, on a tubular neighborhood of $M$ in $X$, the tangent bundle $TX$ can be identified with $E$ over $(-1, 0] \times M$.

It is well-known that the connection $\omega$ (with curvature $\Omega$) in Section 1 can be extended to an $SO(n+1)$-connection, which is still denoted by $\omega$ (with curvature $\Omega$), on $TX$. Also notice that the restriction of the tangent unit sphere bundle of $X$, denoted by $STX$, to M is $\Sigma M$. Let $\bar{\sigma}: STX \to X$ be the projection, which extends $\sigma$.

Now let $V$ be a nowhere zero smooth vector field on $M$ which is tangent to $X$, and let $\overline{V}$ be a generic extension of $V$ on $X$. Without loss of generality, we may assume $\overline{V}$ has only one zero point $p$.

For $r > 0$, let $B_r(p)$ be the geodesic ball of radius $r$ around $p$. Then for small $r$ (when $B_r(p)$ is in the interior of $X$), $\overline{V}$ naturally defines a cross section $\bar{\alpha}: X \setminus B_r(p) \to STX$, which restricts to $\alpha$ on $M$.

Assume first that $n$ is odd; it follows from Proposition 1.1:

$$-\chi(X) = -\int_X E(\Omega) = -\lim_{r \to 0^+} \int_{X \setminus B_r(p)} \bar{\alpha}^* \bar{\sigma}^*(E(\Omega)) = \lim_{r \to 0^+} \int_{X \setminus B_r(p)} d\bar{\alpha}^*(\Psi)$$

$$= \int_M \alpha^*(\Psi) - \lim_{r \to 0^+} \int_{\partial B_r(p)} \bar{\alpha}^*(\Psi) = \int_M \alpha^*(\Psi) - \operatorname{ind}_p(\overline{V})$$

where the first equality follows from the Gauss-Bonnet theorem, the second follows from the fact that $\bar{\sigma}\bar{\alpha} = \operatorname{id}$, and the fourth is by Stokes' theorem.

Theorem 0.1 then clearly follows when $n$ is odd. The case when $n$ is even is similar. If $X$ is not orientable, from the proof above, the theorem easily follows by passing to the orientable double covering of $X$. The proof is therefore complete.

Some special cases worth mentioning are:

- When $V$ is transversal to $M$, it is easy to see $\alpha^*(\Psi) = 0$ if $n$ is odd, while $\alpha^*(\Psi) = \frac{1}{2}$ times the Euler curvature form of $TM$ if $n$ is even (and if $V$ is pointing out of $X$).

- When $V$ is tangent to $M$, it is easy to see $\alpha^*(\Psi) = 0$ for both odd and even $n$.

The corresponding formulae for $\sum \operatorname{ind}_{p_j}(\overline{V})$ in these cases can also be seen easily from the Poincaré-Hopf theorem, except maybe one—when $n$ is even and $V$ is tranversal to $M$, which is the relative Poincaré-Hopf theorem (cf. [P]).

It is interesting to compare our formula with the one in [G] or [P]. This yields
$$\alpha^*(\Upsilon)([M]) = \begin{cases} -\operatorname{Ind}(\partial_- V) & \text{if } n \text{ is odd} \\ \chi(X) - \operatorname{Ind}(\partial_- V) & \text{if } n \text{ is even} \end{cases}.$$

We refer to [G] and [P] for the definition of $\operatorname{Ind}(\partial_- V)$.

## Section 3

We now turn to the general oriented vector bundle $\xi$. Let $\alpha_0 : M \to \hat{\xi}$ be the canonical $\infty$-cross section, and as before $\imath : S^n \to \hat{\xi}$ be any (orientation-preserving) diffeomorphism from $S^n$ into a fibre of $\sigma$.

By Proposition 1.1 and a special case mentioned at the end of Section 2, the element $\Upsilon(\xi) \in H^n(\hat{\xi}, \mathbb{R})$ represented by $\Psi$ constructed in Section 1 has the following properties:

(1) $\imath^*(\Upsilon(\xi)) = s^n$, where $s^n$ denotes the canonical generator of $H^q(S^n, \mathbb{R})$.

(2) $\alpha_0^*(\Upsilon(\xi)) = -\frac{1}{2}e(\xi)$, where $e(\xi) \in H^n(M, \mathbb{R})$ is the real coefficient Euler class of $\xi$.

*Example.* Let $M = S^2$, and let $\xi = TS^2$ and $\eta = M \times \mathbb{R}^2$ be the trivial (oriented) plane bundle over $S^2$. Then topologically $\hat{\xi} = \hat{\eta} = S^2 \times S^2$. Let $i_k : S^2 \times S^2 \to S^2$, $k = 1, 2$, be the projections onto the two factors respectively. It is seen immediately from the construction in Section 1 that $\Upsilon(\xi) = i_1^*(s^2) + i_2^*(s^2)$ and $\Upsilon(\eta) = i_2^*(s^2)$.

Guided by (1), (2) above, we now define $\Upsilon(\xi)$ without using the connections.

PROPOSITION 3.1.  *The following sequence*
$$0 \longrightarrow H^n(M, \mathbb{Z}) \xrightarrow{\sigma^*} H^n(\hat{\xi}, \mathbb{Z}) \xrightarrow{\imath^*} H^n(S^n, \mathbb{Z}) \longrightarrow 0$$

*is exact.*

*Proof.* The proposition comes easily from the following commutative diagram of the Gysin sequence (cf. [MS])

$$\begin{array}{ccccccccc} 0 & \longrightarrow & H^n(M) & \xrightarrow{\sigma^*} & H^n(\hat{\xi}) & \longrightarrow & H^0(M) & \longrightarrow & 0 \\ & & & & \downarrow \imath^* & & \downarrow \approx & & \\ & & & & H^n(S^n) & \xrightarrow{\approx} & H^0(\text{point}) & & \end{array}$$

where the integer coefficients are understood. The first horizontal line, which is exact, is from the Gysin sequence of the vector bundle $\nu \oplus \xi$. As before $\nu$ is the canonical trivial oriented line bundle, and we have used the fact that $e(\nu \oplus \xi) = 0$ to conclude that the homomorphism $H^0(M) \to H^{n+1}(M)$ in the Gysin sequence vanishes. $\square$

Proposition 3.1 easily implies that there is a canonical decomposition

$$H^n(\hat{\xi}, \mathbb{Z}) = \sigma^*(H^n(M, \mathbb{Z})) \oplus \alpha_0^{*-1}(0)$$

and $\imath^*|_{\alpha_0^{*-1}(0)} : \alpha_0^{*-1}(0) \to H^n(S^n, \mathbb{Z})$ is an isomorphism. Needless to say $\alpha_0^*|_{\sigma^*(H^n(M,\mathbb{Z}))} : \sigma^*(H^n(M, \mathbb{Z})) \to H^n(M, \mathbb{Z})$ is also an isomorphism.

We can now define $\Upsilon(\xi) \in H^n(\hat{\xi}, \mathbb{Z}) \otimes \frac{1}{2}$, where $H^n(\hat{\xi}, \mathbb{Z}) \otimes \frac{1}{2}$ denotes the tensor product, as $\mathbb{Z}$-module, of $H^n(\hat{\xi}, \mathbb{Z})$ and the subgroup of $\mathbb{Q}$ generated by $\frac{1}{2}$, as follows:

$$\Upsilon(\xi) = -\frac{1}{2}\sigma^*(e(\xi)) + \imath^*|_{\alpha_0^{*-1}(0)}^{-1}(s^n).$$

Since the sequence in Proposition 3.1 is clearly also exact with real coefficient, properties (1) and (2) above characterize $\Upsilon(\xi)$, defined in Section 1, in $H^n(\hat{\xi}, \mathbb{R})$. Obviously, this agrees with the $\Upsilon(\xi)$ just defined in this section, after tensoring with $\mathbb{R}$. This shows that the element $\Upsilon(\xi) \in H^n(\hat{\xi}, \mathbb{R})$ constructed as in Section 1 does not depend on the choice of connections.

It is well-known that if an oriented $M$ is the boundary of a compact manifold, then $e(TM) \in H^n(M, \mathbb{Z})$ is even. Hence in this case (also in the case $n$ is odd) $\Upsilon \in H^n(\Sigma M, \mathbb{Z})$.

To finish, let us now prove Theorem 0.2 from the introduction. Here again we use the integer coefficients.

First, it follows immediately, from the Gysin sequence of $\nu \oplus \xi$, that $\sigma^* : H^{n-1}(M) \to H^{n-1}(\hat{\xi})$ is an isomorphism. Hence so is $\alpha_0^* : H^{n-1}(\hat{\xi}) \to H^{n-1}(M)$.

Then from the cohomology exact sequence of the pair $(\hat{\xi}, \xi_\infty)$,

$$\cdots \longrightarrow H^{n-1}(\hat{\xi}) \xrightarrow{\alpha_0^*} H^{n-1}(M) \longrightarrow H^n(\hat{\xi}, \xi_\infty) \xrightarrow{\jmath^*} H^n(\hat{\xi}) \xrightarrow{\alpha_0^*} H^n(M) \longrightarrow \cdots$$

where we have replaced $H^j(\xi_\infty), j = n - 1, n$ by $H^j(M)$, we see that $\jmath^* : H^n(\hat{\xi}, \xi_\infty) \to H^n(\hat{\xi})$ is injective, and its image is $\alpha_0^{*-1}(0)$.

By the definition of $\Upsilon(\xi)$, to prove Theorem 0.2, it is now sufficient to verify $\imath^*(\jmath^*(\gamma(\xi)))$ as the canonical generator of $H^n(S^n)$. But this easily follows from the characterization of the Thom class $\gamma(\xi)$.

INDIANA UNIVERSITY, BLOOMINGTON, IN
*E-mail address*: jsha@indiana.edu

## REFERENCES

[C] S. S. CHERN, A simple intrinsic proof of the Gauss-Bonnet formula for closed Riemannian manifolds, Ann. of Math. **45** (1944), 747–752.
[CS] S. S. CHERN and J. SIMONS, Characteristic forms and geometric invariants, Ann. of Math. **99** (1974), 48–69.
[G] D. H. GOTTLIEB, The law of vector fields, preprint.
[KN] S. KOBAYASHI and K. NOMIZU, *Foundations of Differential Geometry* I, II, Interscience, New York, 1963, 1969.
[M] J. MILNOR, *Topology from the Differentiable Viewpoint*, The Univ. of Virginia Press, Charlottesville, VA, 1965.
[MS] J. MILNOR and J. D. STASHEFF, *Characteristic Classes*, Ann. of Math. Studies, No. 76, Princeton University Press, Princeton, NJ, 1974.
[P] C. C. PUGH, A generalized Poincaré index formula, Topology **7** (1968), 217–226.